\documentclass[11pt,onecolumn]{article}
\usepackage{graphicx} % Required for inserting images
\usepackage[T1]{fontenc}
\usepackage{lmodern}
\usepackage{microtype}
\usepackage{authblk}
\usepackage[bottom=2cm]{geometry}
\usepackage{amsmath}
\usepackage{amssymb}
\usepackage{amsfonts}
\usepackage{enumitem}
\usepackage{tikz-cd}
\usepackage{mathtools}
\usepackage{systeme}
\usetikzlibrary{babel}
\usepackage[english]{babel}
\usepackage{amsthm}
\usepackage{amsrefs}
\usepackage{aliascnt}
\usepackage[onehalfspacing]{setspace}
\usepackage{hyperref}
\usepackage{cleveref}
\tikzcdset{scale cd/.style={every label/.append style={scale=#1},
    cells={nodes={scale=#1}}}}
\newcommand{\N}{{\mathbb{N}}}

\newcommand{\C}{{\mathbb{C}}}

\newcounter{ex}

\setlist[enumerate]{label=(\roman*),itemsep=3pt}
\theoremstyle{definition}
\newtheorem{definition}{Definition}[section]
\crefname{definition}{definition}{definitions}
\Crefname{definition}{Definition}{Definitions}
\newaliascnt{exmp}{definition}

\aliascntresetthe{exmp}
\crefname{exmp}{example}{examples}
\Crefname{exmp}{Example}{Examples}
\newaliascnt{theorem}{definition}
\newtheorem{theorem}[theorem]{Theorem}
\aliascntresetthe{theorem}
\crefname{theorem}{theorem}{theorems}
\Crefname{theorem}{Theorem}{Theorems}
\newaliascnt{corollary}{definition}
\newtheorem{corollary}[corollary]{Corollary}
\aliascntresetthe{corollary}
\crefname{corollary}{corollary}{corollaries}
\Crefname{corollary}{Corollary}{Corollaries}
\newaliascnt{remark}{definition}
\newtheorem{remark}[remark]{Remark}
\aliascntresetthe{remark}
\crefname{remark}{remark}{remarks}
\Crefname{remark}{Remark}{Remarks}
\newaliascnt{lemma}{definition}
\newtheorem{lemma}[lemma]{Lemma}
\aliascntresetthe{lemma}
\crefname{lemma}{lemma}{lemmas}
\Crefname{lemma}{Lemma}{Lemmas}
\newaliascnt{prop}{definition}

\aliascntresetthe{prop}
\crefname{prop}{proposition}{propositions}
\Crefname{prop}{Proposition}{Propositions}

\title{Constructive Quantifier Elimination with a Focus on Matrix Rings}
\author{Max Illmer and Tim Netzer}
\affil{Department of Mathematics, University of Innsbruck}
\date{\today}

\begin{document}
\maketitle
{\centering\emph{To the memory of Alexander Prestel}\par}
\begin{abstract}

We give a sufficient condition for a model-theoretic structure $B$ to ``inherit'' quantifier elimination from another structure $A$. This yields an alternative proof of one of the main results from \cite{kle}, namely quantifier elimination for certain matrix rings. The original proof uses model theory, and while it is very elegant and insightful, the proof we propose is much shorter and provides a constructive algorithm.
\end{abstract}

\begin{section}{Introduction and Preliminaries}
\subsection{Introduction}
In model theory it is common to infer certain properties of a structure from another, by giving a strong enough relation between the two. The two structures $k$ and ${\rm Mat}_n(k)$, where $k$ is a field and ${\rm Mat}_n(k)$ denotes the matrix ring, are interesting examples of this.  Although they are bi-interpretable after naming matrix units in the matrix ring, ${\rm Mat}_n(k)$ does in general not inherit quantifier elimination from $k$, even if $k$ is, for example, real closed or algebraically closed. However, extending the language of (ordered) rings by the two new functions trace and conjugate transpose establishes a strong enough relation between the structures for quantifier elimination to actually transfer. This was first proven in  \cite{kle} using model completeness and an amalgamation argument for the expanded matrix structure.

In this paper, we first discuss such relations between structures in broad generality, and then apply them to matrix rings. The  main idea is to find maps that have enough expressive power to transfer a formula into the structure with quantifier elimination, and then pull it back again.  We will use the same algebraic ingredient as in \cite{kle}, the so-called Specht property of certain fields. This provides maps $f_m: {\rm Mat}_n(k)^m \to k^{d(m)}$ from matrices into the finite-dimensional  vector space $k^{d(m)}$, whose fibers determine matrix-tuples up to simultaneous unitary similarity, which is granular enough for any definable set in ${\rm Mat}_n(k)^m$. The only thing to show then is that the subset of $k^{d(m)}$, which ``parametrizes'' the definable set of interest in ${\rm Mat}_n(k)^m,$ is itself definable in the structure $k$. Important to note here is that if quantifier elimination over $k$ can be done in a constructive way (as is the case for example over algebraically closed  fields or real closed fields) and the invariants $f_m$ are given explicitly,  quantifier elimination can also be done constructively for ${\rm Mat}_n(k)$ by our approach.

\subsection{Preliminaries from Algebra}
Since we mostly deal with algebraic structures in this work, we introduce the most important notions.

Throughout, rings are unital and $\N=\{1,2,\ldots\}$. Let $k$ be a field. By ${\rm Mat}_n(k)$ we denote the ring of $n\times n$-matrices with entries from $k,$ which is a non-commutative ring for $n\geqslant 2.$ The usual {\it transpose} of a matrix $A$ is denoted by $A^t$, and its {\it trace}, i.e.\ the sum of its diagonal entries, by ${\rm tr}(A).$

The field $k$ is {\it algebraically closed,} if every non-constant polynomial over $k$ has a zero in $k$. The most prominent example is $k=\mathbb C$, but every field is contained in an algebraically closed field.

The field $k$ is {\it real closed,} if $-1$ is not a square in $k,$ and the field extension $k(\sqrt{-1})$ is algebraically closed (in which case it is the algebraic closure of $k$). There are other characterizations \cite{pd}, but this one is the easiest to write down. Here, $\mathbb R$ is the most prominent example. Real closed fields have a unique field ordering, where an element is nonnegative if and only if it is a square. Every ordered field is contained in a real closed field. 

If $k$ is real closed, then on  $k(\sqrt{-1})$ one obtains an {\it involution}, just as complex conjugation on $\C,$ by multiplying the coefficient of $\sqrt{-1}$ by $-1$. Thus on ${\rm Mat}_n(k(\sqrt{-1}))$ we do not only have the transposition operation, but also {\it conjugate transposition}, denoted by $A^*$ for a matrix $A$, which combines transposition with entry-wise conjugation.  For matrices over $k$ itself, we set $A^*=A^t$. This further gives rise to the notion of {\it self-adjoint} ($A^*=A$) and {\it unitary} ($A^*=A^{-1}$) matrices.

A self-adjoint matrix $A$ is {\it positive semidefinite}, if all its eigenvalues are nonnegative elements (of the real closed field $k$). The usual alternative characterizations of positive-semidefiniteness remain valid for real closed fields $k$ and their algebraic closures $k(\sqrt{-1}).$ 

\subsection{Preliminaries from Logic}
Throughout this paper we will work in the classical setup of first order logic and model theory, as for example described in \cite{cha}, \cite{hod}  and \cite{pre}. We denote formal languages by calligraphic letters such as $\mathcal{L}$.
A language can be extended by an additional symbol $S$, we denote this by  $(\mathcal{L},S)$. 
In order to speak of rings, a language $\mathcal{L}$ needs at least two function symbols $+,\cdot$ of arity $2$, and two constant symbols $0,1$. With these symbols one can write down the ring axioms, and the axioms of algebraically closed or real closed fields, for example. Throughout this paper, we will denote this specific language by $\mathcal{F}$. The complex numbers $\C$ form an obvious $\mathcal{F}$-structure, with the usual interpretation of $+,\cdot, 0,1$.
If we want to talk about ordered rings, we can extend the language by the relation symbol $\leqslant$. An obvious $(\mathcal{F},\leqslant)$-structure is $\mathbb{R}$ with the usual interpretations. Note that for real closed fields, the extension by $\leqslant$ is to some extent redundant, since nonnegative elements are squares, and nonnegativity can thus be expressed in $\mathcal F$ alone. However, this uses a quantifier, and to obtain  quantifier elimination one needs to include the relation symbol $\leqslant$.

To emphasize the free variables $x_1,...,x_m$ of a formula $\phi$, we also sometimes write $\phi(x_1,...,x_m)$. Given such a formula, and a structure $A$, $\phi$ defines the subset of $A^m$  of all $m$-tuples that fulfill the formula, and we write $\phi(A)$ for this subset. Such sets are called {\it definable}. Unless stated otherwise, definability is understood without parameters.

An $\mathcal L$-structure $A$ admits \textit{quantifier elimination}, if every definable set is also definable without quantifiers. More precisely, if  for every $\mathcal L$-formula $\phi(x_1,...,x_m)$ there exists a quantifier-free $\mathcal L$-formula $\psi(x_1,...,x_m)$ with $\phi(A)= \psi(A).$

It is a well-known fact that algebraically closed fields, i.e.\ $\mathcal{F}$-structures satisfying the axioms for algebraically closed fields, and real closed fields, i.e.\ ($\mathcal{F},\leqslant)$-structures satisfying the axioms of real closed fields, admit quantifier elimination. In fact a stronger result holds, namely that quantifiers can be removed independently of the actual field (this means that the \textit{theory} of algebraically or real closed fields admits quantifier elimination).
Both results can actually be proven by showing that the projection of a quantifier-free definable set is again quantifier-free definable (quantifier-free definable sets for algebraically closed fields are constructible sets in the Zariski topology, while for real closed fields they are the semialgebraic sets). Since projections are defined by formulas of the form $\exists x\ \phi(x,x_1,...,x_m)$, the general case can then be deduced inductively. In the case of algebraically closed fields, a constructive proof boils down to a simple application of the Euclidean algorithm. For real closed fields the result is less trivial and is the content of the Tarski-Seidenberg Theorem. But also here there exist constructive proofs, see for example \cites{basu, boch, sch}.

\subsection{Complete Invariants}
\begin{definition}
Given an equivalence relation $\sim$ on some set $X$, we say that a function $f\colon X\to Y$ is a \textit{complete invariant}, if the following holds true: $$\forall  a,b\in X\colon\quad  a \sim b \iff f(a)=f(b).$$
It is also common to talk of a \textit{complete set of invariants}, if $f$ and $Y$ can be written as a direct product in a compatible way.
\end{definition} 

Very often, additional properties for the set $Y$ and the function $f$ are required in order to make good use of the invariants. In our proof below for example, we want $Y$ to be a finite-dimensional vector space over a field $k$.

One crucial ingredient in the proof from \cite{kle} is the following theorem, which characterizes simultaneous unitary similarity of matrices  in terms of  a very nice complete set of invariants.
\begin{theorem}\label{Spe}
    Let $k$ be a real closed field or its algebraic closure. Then for  $n,m\in \N$ and $A=(A_1,\ldots, A_m),B=(B_1,\ldots, B_m)\in {\rm Mat}_n(k)^m$ 
    there exists a unitary $U\in {\rm U}_n(k)$ with $$ U^*A_iU=B_i \ \mbox{ for } \ i=1,\ldots,m$$ if and only if $${\rm  tr}(\omega(A_1,...,A_m,A_1^*,...,A_m^*))={\rm tr}(\omega(B_1,...,B_m,B_1^*,...,B_m^*))$$
holds for all words $\omega$ of length $\leqslant n^2$ in $2m$ non-commuting variables.
\end{theorem}
For classical results on simultaneous unitary similarity and matrix invariants, see \cites{wie,sib}. The finite word-length bound used here follows from the invariant-theoretic results in \cite{pro}*{Theorem 7.3} and \cite{raz}. For each fixed $n,m$, the resulting statement over $\mathbb R$ or $\mathbb C$ transfers to arbitrary real closed fields and their algebraic closures by expressing matrix entries in real coordinates and using completeness of the theory of real closed fields.

If such a bound depending on $n$ and $m$ exists for a field $k$, it is also said that $k$ has the \textit{Specht property}. So the above theorem says that real closed fields and their algebraic closures have the Specht property.

Why are we interested in the specific equivalence relation of simultaneous unitary similarity? The reason is its natural connection with the language $(\mathcal{F},{\rm tr},{}^*,\leqslant)$.
If $k$ is a real closed field or its algebraic closure, we consider ${\rm Mat}_n(k)$ as an $(\mathcal{F},{\rm tr},{}^*, \leqslant)$-structure, where  ${\rm tr}$ is interpreted as the usual trace (multiplied with the identity),  ${}^*$ as conjugate transposition, and $\leqslant$ as the partial ordering $A\leqslant B$ if and only if $B-A$ is positive semidefinite.

Note that for all results to come, interpreting $\leqslant$ as the order on the (self-adjoint) center would also suffice. In fact together with trace and adjoint it can be used to define positive semidefiniteness without quantifiers (see for example below, or also \cite{cle}). We instead opt for the  interpretation of $\leqslant$ by positive semidefiniteness, which is more natural in the context of free semialgebraic geometry.

\begin{lemma}\label{Lem}
For the $(\mathcal{F},{\rm tr},{}^*,\leqslant)$-structure ${\rm Mat}_n(k)$,  definable sets are closed under simultaneous unitary similarity.
\end{lemma}
\begin{proof}For any $U\in {\rm U}_n(k)$ and $A,B\in {\rm Mat}_n(k)$ we have \begin{align*}U^*AU+U^*BU&=U^*(A+B)U\\
 U^*AUU^*BU&=U^*ABU\\
 {\rm tr}(U^*AU)&={\rm tr}(A)=U^*{\rm tr}(A)U\\
 (U^*AU)^*&=U^*A^*U\end{align*} and
 \begin{align*}
A=0 &\iff U^*AU=0\\ A\geqslant 0 &\iff U^*AU\geqslant 0.\end{align*}
If  $D\subseteq{\rm Mat}_n(k)^m$ is defined by a quantifier-free formula $\phi(x_1,...,x_m)$, this  easily implies  $(A_1,...,A_m)\in D \Longrightarrow (U^*A_1U,...,U^*A_mU)\in D $. Existential quantifiers  correspond to projections, and this clearly preserves closedness under simultaneous unitary similarity. For general formulas the statement then follows by induction. 
\end{proof}
\begin{remark}
For the readers with a deeper background in model theory, this is just the trivial fact that the maps $\phi_U\colon {\rm Mat}_n(k)\to {\rm Mat}_n(k), A \mapsto U^*AU$ for $U \in U_n(k)$ form a subgroup of  the automorphism group ${\rm Aut}({\rm Mat}_n(k)),$ and thus all elements in the orbits $\{\phi^m_U(X) , U\in U_n(k)\}$ are of the same complete type over $\emptyset$ (indeed, even over $k\cdot I_n$, which these automorphisms fix pointwise).

By contrast, the maps $\phi_G(A)=G^{-1}AG$ for $G\in {\rm GL}_n(k)$ preserve the ring and trace operations, but generally not adjoint or the positive-semidefinite order. For example, conjugating $\begin{psmallmatrix}1&1\\1&1\end{psmallmatrix}$ by $\operatorname{diag}(1,2)$ gives $\begin{psmallmatrix}1&2\\1/2&1\end{psmallmatrix}$, which is not self-adjoint. Thus simultaneous similarity is not an automorphism-orbit relation for the full structure considered here. Ranks of linear matrix pencils do separate simultaneous similarity orbits \cite{ranks}, but this alone does not verify the quantifier-free pullback condition ($ii$) of \Cref{theo}. In particular, traces of powers do not determine rank: a nonzero nilpotent matrix has the same power traces as zero. Using $A^*A$ instead recovers rank from its spectrum, but uses the adjoint symbol. The distinctions between the ring-and-trace reduct, ordering on the center, and the full positive-semidefinite order must therefore be kept in mind when comparing with \cite{kle}.
\end{remark}
\end{section}

\begin{section}{Results}

We start by describing in a rather general way a relation between two structures, which allows to transfer quantifier elimination. 
To this end, let $\mathcal{L}$ and $\mathcal{G}$ be formal languages and $A$ and $B$ structures over $\mathcal{L}$ and $\mathcal{G},$ respectively. We further assume that for every $m\in \N$ there exists a function $f_m\colon B^m \to A^{d(m)}$ with the following three properties:
\begin{enumerate}
    \item[($i$)] For any definable set $D\subseteq B^m$, the set $f_m(D)$ is definable in $A^{d(m)}.$
    \item[($ii$)] For any  quantifier-free definable set  $E\subseteq A^{d(m)}$, the set $f_m^{-1}(E)\subseteq B^m$ is quantifier-free definable.
    \item[($iii$)] For any $a\in A^{d(m)}$ and any definable set $D\subseteq B^m$, we either have $f_m^{-1}(a)\subseteq D$ or $f_m^{-1}(a)\cap D=\emptyset.$
\end{enumerate}
We call conditions ($i$) and ($ii$) {\it constructive,} if there is a constructive algorithm that turns  definitions of sets into definitions of their image/preimage, as required.

Properties ($i$) and ($ii$) allow to transform definable sets back and forth via $f_m$. Property ($iii$) basically says that $f_m $ is ``injective up to definability''.

\begin{theorem}\label{theo}
   Let $A,B$ and $(f_m)_{m\in \N}$ be as described above. Then, if $A$ admits quantifier elimination, so does $B$.  Furthermore, if  ($i$) and ($ii$) are constructive, and there exists a constructive algorithm for quantifier elimination in $A$, then there is also a constructive algorithm for quantifier elimination in $B$.
\end{theorem}
\begin{proof}
Let $D\subseteq B^m$ be definable. By ($i$) we know that $f_m(D)$ is definable in $A^{d(m)}$. Since $A$ admits quantifier elimination, this set is also quantifier-free definable, and therefore by ($ii$) so is $f_m^{-1}(f_m(D))$. But due to ($iii$) we have $f_m^{-1}(f_m(D))=D$, which finishes the proof. The statement about constructiveness follows directly from this argument.
\end{proof}

Before getting to our main application, let us demonstrate the use of \Cref{theo} in some easy cases.

\begin{corollary}\label{cor:conj}
 Let $k$ be a real closed field and $K=k(\sqrt{-1})$ its algebraic closure. Then $K$ admits constructive quantifier elimination in $(\mathcal F,\bar{\cdot},\leqslant)$, where $\bar{\cdot}$ is conjugation and $\leqslant$ is the order on $k$.
\end{corollary}
\begin{proof}
 For $z\in K^m$, consider
 \[
 P_z(t,x_1,\ldots,x_m)=
 \left(t+\sum_{j=1}^m x_jz_j\right)
 \left(t+\sum_{j=1}^m x_j\overline{z}_j\right).
 \]
 Let $f_m(z)$ be its coefficient tuple, omitting the coefficient of $t^2$. Its entries are $z_j+\overline{z}_j$, $z_j\overline{z}_j$, and $z_j\overline{z}_\ell+\overline{z}_jz_\ell$ for $j<\ell$. All are fixed by conjugation, so $f_m$ takes values in $k^{m(m+3)/2}$.

 We verify the hypotheses of \Cref{theo}. For ($i$), write elements of $K$ as $x+y\sqrt{-1}$ with $x,y\in k$. Field operations, conjugation, and the restricted order have polynomial coordinate descriptions over $k$. Images under $f_m$ are then definable by existentially quantifying these coordinates. Since every entry of $f_m$ is a term in the expanded field language with values in $k$, ($ii$) follows by direct substitution.

 Finally, if $f_m(z)=f_m(w)$, then $P_z=P_w$. Unique factorization in $K[t,x_1,\ldots,x_m]$ forces the monic linear factor $t+\sum_jx_jw_j$ to equal one of the two factors of $P_z$. Comparing coefficients gives $w=z$ or $w=\overline{z}$, with one common choice for the entire tuple. Conversely, simultaneous conjugation leaves $P_z$ unchanged. Thus the fibers are precisely the simultaneous-conjugation orbits. These cannot be split by definable sets, since conjugation preserves the full structure. This proves ($iii$). The conclusion, including constructiveness, follows from \Cref{theo} and constructive quantifier elimination over the real closed field $k$.
\end{proof}

\begin{remark}\label{conj}
The preceding proof (as well as the coming one for \Cref{Mat}) already illustrates the invariant-theoretic approach for the transfer functions, further discussed in \cite{illmer-invariants}*{Propositions 2.7 and 2.9}. The same argument gives quantifier elimination for quadratic Galois extensions $K=k(a)$ in arbitrary characteristic if $k$ admits quantifier elimination in $\mathcal L\supseteq\mathcal F$ and the minimal-polynomial coefficients of $a$ are definable in $k$. The defining coefficients may be named in $\mathcal L$ so no constant naming $a$ is needed. More generally, for finite Galois $K/k$ with group $G$, the coefficients of $\prod_{g\in G}(t+\sum_jx_jg(z_j))$ separate simultaneous $G$-orbits. Thus \Cref{theo} applies provided coordinate descriptions are definable over $k$, $k$ admits quantifier elimination, $G$ preserves the expanded language, and quantifier-free pullbacks are available. Naming $G$ supplies terms and preserves its action when $G$ is abelian; otherwise only $Z(G)$ preserves all named automorphisms. One may instead name the invariant coefficient functions themselves.
\end{remark}

The following example can also be found as a  preliminary result in \cite{kle}*{Proposition 2.1.11}.

\begin{corollary}
 Let $\mathcal L$ be a language extending  the language $\mathcal F$ of rings, and let $k$ be a ring that admits quantifier elimination w.r.t.\ $\mathcal L$. Consider the ring ${\rm Mat}_n(k)$ as an $(\mathcal{L},(e_{ij})_{i,j \leqslant n})$-structure, where all symbols from $\mathcal L$ beyond $\mathcal F$ are interpreted only on the center and as in $k$ (functions are set $0$ outside), and $e_{ij}$ is interpreted as the usual matrix unit $E_{ij}$. Then ${\rm Mat}_n(k)$   admits quantifier elimination.
\end{corollary}
\begin{proof}
We use the functions  $f_m\colon {\rm Mat}_n(k)^m \to k^{n^2m}$ that map a matrix tuple to the tuple of all its matrix entries. Then property ($i$) is easy to see, ($ii$) follows from the fact that $E_{1i}AE_{j1}=a_{ij}E_{11}$ holds for all $A=(a_{ij})_{i,j}\in{\rm Mat}_n(k)$, and 
 ($iii$) is clear from injectivity of $f_m$.
\end{proof}
This result is intuitively clear, as the addition of matrix units allows us to freely pick any particular matrix entry and operate on it. Much less obvious is the following result, shown in a slightly different form in \cite{kle}.

\begin{corollary}\label{Mat}
Let $k$ be a real closed field or the algebraic closure thereof. Then the $(\mathcal{F},{\rm tr},{}^*,\leqslant)$-structure ${\rm Mat}_n(k)$ admits quantifier elimination.
\end{corollary}
\begin{proof}
Consider the functions $f_m\colon {\rm Mat}_n(k)^m \to k^{d(m)}$ that  compute the invariants from the Specht property, explicitly $$f_m(A_1,...,A_m)=({\rm tr}(\omega(A_1,...,A_m,A_1^*,...,A_m^*)))_{\omega \in J(n,m)}$$ where $J(n,m)$ is the set of words of length $\leqslant n^2$ in $2m$ variables, including the empty word (evaluated as $I_n$). Here the coordinates of $f_m$ are ordinary scalar traces; their scalar-matrix versions are terms in the matrix language.
Since $k$ admits quantifier elimination as an $(\mathcal{F},{}^*,\leqslant)$-structure with order restricted to the fixed field of ${}^*$ (see \Cref{cor:conj} for the algebraically closed case), it only remains to show that $f_m$ fulfills properties ($i$), ($ii$) and ($iii$), For ($i$), translate matrices into entries; positive semidefiniteness is expressed by self-adjointness and nonnegative principal minors. This gives definable coordinate descriptions, and images under $f_m$ are obtained by existentially quantifying the entries. For ($ii$) let $E\subseteq k^{d(m)}$ be defined by the quantifier-free formula $\varphi$ in the  free variables $x_\omega$ (for $\omega\in J(n,m)$).
Substitute ${\rm tr}(\omega(A_1,...,A_m,A_1^*,...,A_m^*))$ for $x_\omega$, translating each ordering atom $S\leqslant T$ between the resulting scalar terms as $S=S^*\land T=T^*\land S\leqslant T$. This gives a quantifier-free formula defining $f_m^{-1}(E).$ ($iii$) is due to the fact that $f_m$ is a complete invariant of simultaneous unitary similarity (see \Cref{Spe}), and the orbits of simultaneous unitary similarity can never be split by a $(\mathcal{F},{\rm tr},{}^*,\leqslant)$-definable set, as shown in \Cref{Lem}.
\end{proof}

\begin{remark}
 We have made a conscious effort to  state and prove the above results in as elementary a way as possible, and indeed for the general statement \Cref{theo} we believe this to be the most natural choice. But from the viewpoint of logic and model theory, the following perspective on our proof of \Cref{Mat} could be of some value:
 We make use of the definable equivalence relation: \begin{align*}(&A_1,...,A_m)\sim_{{\rm U},m} (B_1,...,B_m) \\ & \Leftrightarrow\exists U\colon UU^*=U^*U=I_n \ \land \  UA_1U^*=B_1 \ \land ... \land UA_mU^*=B_m.\end{align*}
 The logical interpretation of \Cref{Spe} then reads as follows: For each $m\in \N$ the {\it imaginaries} $a/\sim_{{\rm U},m}$ can be eliminated uniformly, a unique code being given by $\iota^{d(m)}(f_m(a))$, where $\iota$ denotes the $(\mathcal{F},{}^*)$-embedding of $k$ in ${\rm Mat}_n(k)$ ($\iota(x)=x \cdot I_n$).  Thus the quotients  ${\rm Mat}_n(k)^m/\sim_{{\rm U},m}$ are in bijection to  $\iota^{d(m)}(f_m({\rm Mat}_n(k)^m)) \subset (k\cdot I_n)^{d(m)}$. The coding map is a tuple of terms in the matrix language. Every definable subset of ${\rm Mat}_n(k)^m$ is a union of its fibers, and its image is definable in the scalar-field structure by condition ($i$). Condition ($ii$) then pulls a quantifier-free definition of this image back to the matrix structure.

In the general setup of \Cref{theo}, no joint definability of $f_m$ is assumed, so its fiber relation need not be definable in $B$. A useful sufficient approach is to find a complete invariant for the orbits of a subgroup of ${\rm Aut}(B)$, with values in a structure admitting quantifier elimination, and then check image definability and quantifier-free pullbacks separately. Invariance alone does not suffice: the fibers must not join tuples separated by a definable set.
\end{remark}

In the following, $J(n,m)$ will denote the index set for the co-domain of the invariants $f_m$ from the proof of \Cref{Mat}. Explicitly, we have $$d(m)=|J(n,m)|=\sum_{i=0}^{n^2}(2m)^i=\frac{(2m)^{n^2+1}-1}{2m-1}.$$
Let us now elucidate how the above proof yields an algorithm for quantifier elimination in ${\rm Mat}_n(k)$. 
We are given a formula $\phi(X_1,\ldots, X_m)$ in the language $(\mathcal F,{\rm tr},{}^*,\leqslant)$. We first replace each matrix variable in $\phi$ (also the bound ones) by $n^2$ variables for its entries, and spell out all computations to obtain an $(\mathcal F,{}^*,\leqslant )$-formula $\tilde\phi$ in $mn^2$ free commuting variables. If the entries of $X_i$ are denoted $x^i_{rs}$, the free variables of $\tilde\phi$ are all $x^i_{rs}$.  We then consider the formula $$\exists x^1_{rs},\ldots, x^{m}_{rs}\colon \tilde\phi(x^i_{rs}) \wedge x=f_m(x^i_{rs}),$$ where $x$ is a $d(m)$-tuple and $f_m(x^i_{rs})$ denotes the function values computed in terms of the matrix entries, which can indeed be written down in the language $(\mathcal F,{}^*,\leqslant)$. This formula now defines the image of the set defined by $\phi$ under $f_m$.  By quantifier elimination in $k$, there exists an equivalent  quantifier-free $(\mathcal F,{}^*,\leqslant)$-formula $\gamma(x)$. Replacing $x$ by $f_m(X_1,\ldots, X_m)$ and translating ordering atoms as above yields a quantifier-free $(\mathcal F,{\rm tr},{}^*,\leqslant)$-formula $\psi$ which is equivalent to $\phi.$

Note that  the step of calculating and eliminating quantifiers of the ``parameterizing'' set $f_m(\phi({\rm Mat}_n(k))) \subseteq k^{d(m)}$ has been done simultaneously here. 
For the complexity estimate, assume that $k$ is real closed and that $\phi=\exists Y_1,\ldots,Y_l\colon\xi$, where $\xi$ is a conjunction of $r$ atomic formulas with $m$ free matrix variables and integer coefficients. Suppose the matrix terms have degree at most $s$; trace and transpose do not increase entry degrees. Put
\[
N=(l+m)n^2,\qquad D=d(m).
\]
The image formula eliminates $N$ scalar variables with $D$ free invariant coordinates. If all atoms are equalities, it uses at most $rn^2+D$ polynomials of degree at most $\max(s,n^2)$. For inequalities, express positive semidefiniteness by symmetry and nonnegativity of all principal minors. Each matrix atom then contributes at most $n^2+2^n-1$ scalar conditions, of degree at most $ns$. Thus safe general bounds are
\[
S=r(n^2+2^n-1)+D,\qquad \delta=\max(2,ns,n^2).
\]
Single-block quantifier elimination gives an arithmetic complexity bound of the form
\[
(S\delta)^{\mathcal O((N+1)(D+1))}
\]
for this scalar elimination step, by \cite{basu}*{Chapter 14}. This is an arithmetic-operation bound, not a bit-complexity estimate or a lower bound. Since $D\leqslant 2(2m)^{n^2}$, it is singly exponential in a polynomial in $m$ when $n,l,r,s$ are fixed. The large invariant space makes a direct implementation potentially expensive, but does not establish practical impossibility: the special form of the image equations and redundancy among the invariants may allow substantial reductions. The algebraically closed case requires the additional real-coordinate reduction of \Cref{cor:conj}; the counts above refer only to real closed $k$.

The symbols ${\rm tr}, {}^*$ and $\leqslant$ are  very natural for (free) semialgebraic geometry. For example, they allow to use Newton sums (or the Specht property itself) to get information about the spectrum of a matrix. This was used in \cite{cle} for some explicit quantifier-free description of  interesting sets, e.g.  invertible matrices or positive semidefinite matrices (using only the ordering on the center). 
More generally, properties in free semialgebraic geometry such as unitarity, positivity, having a certain rank and the like, are definable  and can thus be described  constructively without  quantifiers, by the construction from above. For practical computations, it is therefore useful to exploit the structure of the particular problem rather than apply the general construction without simplification. Also note that even in specific cases, where we can simplify the problem enough to make a solution with algebra systems possible, the resulting description from our method will not be the simplest possible or most intuitive one in general.

For example, let us define the real positive semidefinite  $2\times 2$-matrices with our method. As $\leqslant$ is already part of our language, the most obvious defining formula would of course be $X \geqslant 0$. What happens, however, if we strictly follow our algorithm? A well-known quantifier-free definition for positive semidefiniteness in the matrix entries $X=\begin{psmallmatrix}x&y\\y&z\end{psmallmatrix}$ is $$x+z\geqslant 0 \land xz-y^2\geqslant 0.$$
Even in this very simple case, the invariant $f_1$ already has codomain $\mathbb{R}^{31}$. However, due to the predetermined symmetry of $X$, we only really have to consider $4$ dimensions. Using this and proceeding as in our proof, with the help of a computer algebra system like Mathematica gives a final formula:
\begin{align*} {\rm tr}(X)&\geqslant 0, {\rm tr}(X)^2-{\rm tr}(X^2)\geqslant 0, \\{\rm tr}(X)^2-2{\rm tr}(X^2)&\leqslant 0, {\rm tr}(X^2)\geqslant 0,\\
{\rm tr}(X)^3-3{\rm tr}(X){\rm tr}(X^2)+2{\rm tr}(X^3)&=0, {\rm tr}(X)^4-2{\rm tr}(X)^2{\rm tr}(X^2)-{\rm tr}(X^2)^2+2{\rm tr}(X^4)=0, \\ \bigwedge_{\omega \in J(2,1)} {\rm tr}(\omega(X,X^*))&={\rm tr}(\omega(X,X)). \end{align*}
Among symmetric matrices, the first line characterizes positive semidefiniteness. The remaining lines can be replaced by $X=X^*$: the last line includes $\operatorname{tr}(XX^t)=\operatorname{tr}(X^2)$, whose difference is $(x_{12}-x_{21})^2$, and the other identities and inequalities hold automatically for symmetric $2\times2$-matrices.

Redundancy is even clearer when describing the center $Z({\rm Mat}_n(k))$ by eliminating quantifiers in the defining formula $$\forall Y\colon  XY=YX.$$ Strictly following our proof would result in a definable set in $k^{2^{n^2+1}-1}$ as our ``parameterizing'' set. This set is actually in bijection to $k$ itself, as the orbits under similarity obviously collapse on the center.
However, a much easier quantifier-free description in this case is  of course  $nX={\rm tr}(X).$

Let us also quickly test our method for the notorious ``wild'' problem of simultaneous similarity of a pair of $n\times n$-matrices over a real closed field $k$:
$$\exists C,D\colon CD=DC=I \land CA_1D=B_1 \land CA_2D=B_2.$$
Allowing the matrix size to increase, this classification problem contains simultaneous similarity of tuples of arbitrary size. Comparing with our complexity analysis from before, even in the case $n=2$ this results in the $\exists$-block elimination of $24$ variables (from a total of $4705$ variables) from $4697$ polynomial equalities with integer coefficients and maximal degree $4$, illustrating the size of this unreduced invariant-based construction. These counts do not imply that deciding simultaneous similarity of given tuples is intrinsically this expensive.

Finally, note that while all our  results hold for any  dimension $n$, the quantifier-free formulas we construct depend on $n$. Quantifier elimination in the dimension-free setup, which is often considered in free semialgebraic geometry, does not hold in general (see \cites{dre,kle} for example).
\end{section}

\section*{Acknowledgments}
We thank Marcus Tressl for useful comments on a preliminary version of the paper. 

\end{document}